\documentclass[twoside]{article}
\usepackage[dvips]{graphics}

\begin{document}

\centerline\title\textbf{Local and Global Hartogs-Bochner Phenomenon in Tubes}

\author{  Joaquim Tavares}

\footnote{Friday, august 6, 2010  17:14}

\begin{abstract}
A generalization of the Hartogs theorem is proved for
 a class of Tube structures $(\mathrm{M}, \textsc{g},\mathcal{V})$. We assume
 that the intervening commutative Lie algebra $\textsc{g}$ admits
 at least $\mathrm{codim\, }\mathcal{V}$ globally solvable generators.
 We give necessary and sufficient conditions for triviality of the first cohomological group with compact support associated to the Tube structure
 to be trivial. A such global result was previously obtained only when $\mathrm{M}=\textbf{R}^{n}\times\textbf{R}^{m}$ with $\partial/\partial x_{j}$ for
 $j=1,...,m$ generating a Lie subalgebra of $\textsc{g}$.
\end{abstract}

\maketitle

{\bf Introduction.}

\medskip

\baselineskip 20pt

 We start recalling the so called Bochner's extension theorem ([Bo1,2]). It states that if $u$ is a holomorphic function defined in
an open connected set $ \textbf R^{m}+i\,\Omega \subset\textbf
C^{n}$ then it extends as a holomorphic function to the linear
convex envelope $\textbf R^{m}+i\,\widehat\Omega$ of  $\textbf
R^{m}+i\,\Omega$  (one year before that Stein ([St]) proved this
result for $n=2$).
 A kind of local version of the Bochner
extension theorem is found in
  Komatsu
[Ko] where $\textbf R^{m}$ is replaced
 by a ball $\mathrm{B}_{R}$
centered at the origin with radius $R$.
 Later Andronikof ([An]) precise the dependence
between $R$ and the size domain of the extension, namely:

 {{\it Let $\Omega\subset \textbf{R}^{n}$ be a convex
bounded set of dimension $>2$; if $$R-\rho >
\sqrt{2}\, diameter \,(\Omega)$$ then each function holomorphic on a
neighborhood of the tube $\mathrm{B}_{R}\times\partial\Omega $
 has a unique holomorphic extension to a neighborhood of the tube
  $\mathrm{B}_{\rho}\times\Omega $.}}
An example of Ye ([Ye]) shows that $R$ is necessarily bigger than $(1/2)diameter \,(\Omega)$,  leaving the question of finding the sharp constant in the interval $(1/2,\sqrt{2}]$. Another classical
 extension theorem is due to Hartogs([Har1,2]) and it asserts that
a holomorphic function in $\textbf C^{n}\setminus \Omega $, where
$\Omega$ is a bounded open domain with connected boundary
$\partial\Omega$ extends itself to all of $\textbf C^{n}$ as a
holomorphic function. The Bochner extension theorem implies the
Hartogs's one as we see now; let $\textbf{C}^{n}\buildrel
{\Pi}\over \longrightarrow \textbf{R}^{n}$ be the projection into
the imaginary part; $\Pi(x+it)=t$. Then
$$\textbf{C}^{n}\setminus\overline{\Omega}\supset \textbf{R}^{n}+i\,
\textbf{R}^{n}\setminus\Pi(\overline{\Omega})\eqno (0.1)$$ Since
the convex envelope of
$\textbf{R}^{n}\setminus\Pi(\overline{\Omega})$ is
$\textbf{R}^{n}$ the Hartogs extension theorem follows for pairs
$(\Omega, K)$ with $\Omega\setminus{K}$ connected.
Only four years after Fichera([Fi])
published his work reducing the amount of CR data to $\partial\Omega$ under
certain regularity constrains,
 Ehrenpreis ([Eh]) gave a new proof of
the Hartogs extension theorem. The proof of Ehrenpreis was remarkable simple
and its main idea is a cohomological  vanishing argument. The
same idea applied by Hounie \&
Tavares ([HT]) to gives necessary and sufficient conditions for the validity of  the Fichera's version of the
Hartogs extension theorem for a smooth globally integrable
 Tubes structures in $\textbf{R}^{m}\times\textbf{R}^{n}$.
By a smooth globally integrable
Hypoanalytic Tubes structures in $\textbf{R}^{m}\times\textbf{R}^{n}$, we mean a subbundle $\mathcal{L}\subset\textbf{C}\otimes\mathrm{T\,}
(\textbf{R}^{m}\times\textbf{R}^{n})$
such that $\mathcal{L}_{p}=\mathrm{Ker\,}\mathrm{d\, Z}(p)$ where
$\mathrm{Z}:\textbf{R}^{m}\times\textbf{R}^{n}\longrightarrow\textbf{C}^{m}$ is a smooth function function $\mathrm{Z}(x,t)=x+\mathrm{i\,}\Phi(t)$. It extends the concept of the Cauchy-Riemann system in $\textbf{C}^{m}$.
 By a
hypoanalytic structure we understand as a pair $(\mathrm{M}, \mathcal{L})$ consisting of
a smooth manifold $\mathrm{M}$ and a subundle $\mathcal{L}\subset\textbf{C}\otimes\mathrm{T\,M}$ endowed with a associated hypoanalytic atlas $(U_{\alpha}, \mathrm{Z}_{\alpha})$. We mean $\cup_{\alpha}U_{\alpha}=\mathrm{M}$ and the maps  $$\mathrm{Z}_{\alpha}:U_{\alpha}\longrightarrow\textbf{C}
^{m}\,\,\hbox{with}\,\,m=
{\mathrm{dim\,M}-\mathrm{dim}_{\textbf{C}}\mathcal{L}}\eqno(0.2)$$ are smooth and
 $\mathrm{det\,}\mathrm{d\,}\mathrm{Z}_{\alpha}\neq 0$
and if $p\in U_{\alpha}$ then $\mathcal{L}_{p}=\mathrm{ker\,Z}_{\alpha}(p)$.
Finally the etymology comes from the constrain that  $\mathrm{Z}_{\beta}=\mathrm{Z}_{\alpha}\circ \mathrm{H}_{\alpha, \beta}$ in an
neighborhood every point  $p\in U_{\alpha}\cap U_{\beta}$ where $\mathrm{H}_{\alpha, \beta}$ is a biholomorphism
in some open neighborhood of $\mathrm{Z}_{\beta}(p)$.
It is well known that
 $fibers$ of the hypoanalytic structure defined by the germs
$$\mathcal{F}(p)=\mathcal{C}_{\mathrm{Z}_{\alpha}}(p)=
\{\mathrm{Z}_{\alpha}=
\mathrm{Z}_{\alpha}(p)\}\eqno(0.3)$$ are hypoanalytic invariants of the structure. The Sussmann's orbit $\mathcal{O_{\mathcal{L}}}(p)$ (named after Sussmann
([Su])) is the minimal smooth submanifold containing $p$ which supports $\mathcal{L}$ in its complexified tangent space.
We say that a smooth germ of function $u$ at $p$ is $hypoanalytic$ if
$\mathrm{d\,}u$ is a germ of a section of $\mathcal{L}^{\perp}$. If $\mathcal{O_{\mathcal{L}}}(p)$ is compact then the trace of a hipoanalytic function in the orbit must be constant otherwise

             An  Tube structure  $(\mathrm{M}, \mathcal{L}, \mathcal{G})$  is a hypoanalytic structure  endowed
             with a commutative Lie algebra $\mathcal{G}\subset\mathrm{T\,M}$ which verifies the conditions:

             $\bullet_{1}$ if $\mathcal{A}_{p}\subset\mathrm{T}_{p} \mathrm{M}$ is the spam of $\mathcal{G}_{p}$ then $\mathrm{dim\,}\mathcal{A}_{p}\geq\mathrm{codim\,}\mathcal{L}$
            for all $p\in\mathrm{M}$,

             $\bullet_{2}$ $\mathcal{L}_{p}+
             \textbf{C}\otimes\mathcal{G}_{p}=\textbf{C}\otimes
             \mathrm{T}_{p}\mathrm{M}\,\,\,\,\hbox{for all} \,\,p\in \mathrm{M}$,

              $\bullet_{3}$ $[\mathcal{L}, \mathcal{G}]\subset\mathcal{L}$.

             It follows from  $\bullet_{1}$ that $m=\mathrm{dim\,}\mathcal{G}$ is well defined and  greater or equal to $$\mathrm{codim\,}_{\textbf{C}}\mathcal{L}=
             {\mathrm{dim\,M}-\mathrm{dim}_{\textbf{C}}\mathcal{L}}.$$
             Under these hypothesis one can always find an hypoanalytic atlas $(U_{\alpha}, \mathrm{Z}_{\alpha})$ such that
             $\mathrm{Z}_{\alpha}(x,t)=x+\Phi(t)$ for suitable coordinates
             where $\{\partial/\partial x_{1},...,\partial/\partial x_{m}\}$ is a subset of generators of $\textsc{a}$ over $U_{\alpha}$.

             Let us denote by  $\mathcal{F}_{\mathrm{Z}}(p)$ the germ of the closed set $\{\mathrm{Z}=\mathrm{Z}(p)\}$ for an arbitrary hypoanalytic function $\mathrm{Z}$ at $p\in\mathrm{M}$.
 For arbitrary Tubes structures $(\mathrm{M}, \mathcal{L}, \mathcal{G})$ we  have the following characterization of the $local$ Hartogs property;

$\mathbf {Theorem. A.}$ A Tube  $(\mathrm{M}, \mathcal{L}, \mathcal{G})$ has the $local$
             Hartogs property if and only if
             $\mathcal{F}_{\mathrm{Z}}(p)$ is connected for all hypoanalytic germs $\mathrm{Z}$ at $p$ for all
             $p\in \mathrm{M}$.

 {\bf Remark.} Recently was established  in
 work of Henkin \& Michel ([HM]) for abstract real analytic (CR)-structures $(\mathrm{M},\mathcal{L})$ that the  $local$
Hartogs phenomenon is equivalent to
$(\mathrm{M},\mathcal{L})$ be nowhere strictly
pseudoconvex with $\mathrm{dim\,}\mathrm{M}\geq 3$. Actually the concept of pseudoconvexity is belongs to the larger class
of structures called hypoanalytic structures.
The Levi form $\Xi_{\theta}(p)$ is a hypoanalytic invariant defined
in
$\mathcal{L}_{p}\times\mathcal{L}_{p}$ for every
$\theta\in\Sigma_{p}=
\mathcal{L}^{\perp}\cap\mathrm{T}^{\ast}_{p}\mathrm{M}$
by $$\Xi_{p}^{\theta}(\mathrm{v},\mathrm{w})=
              \,\theta \,\big(
              [\mathrm{Re\,}L_{0},\mathrm{Im\,}{
              L_{1}}]\big)(p).$$
             Here   $L_{0},L_{1}$ are germs sections
              of $\mathcal{V}$
              satisfying $L_{0}(p)=\mathrm{v}$ and
              $L_{1}(p)=\mathrm{w}$ and $\Sigma_{p}$ is the characteristic set
              of $(\mathrm{M}, \mathcal{L})$.  Actually it is an well defined object of
              the Sussmman orbit $\mathcal{O}_{\mathcal{L}}(p)$ of
              $\mathcal{V}$. We say that
              $(\mathrm{M},\mathcal{L})$ is strictly
pseudoconvex at $p\in\mathrm{M}$ if
$\Xi_{p}^{\theta}$ is non degenerated with all the eigenvalues with a same sign.   Consequently
               $\mathcal{L}_{p}\subset\textbf{C}
               \otimes\mathrm{T}_{p}
              \mathcal{O}_{\mathcal{L}}(p)$ and
             when $\Xi_{\theta}(p)$ is nondegenerated with all eigenvalues of a same sign we say that $ \mathcal{L}$ is strictly pseudoconvex
             at $\theta\in \Sigma_{p}$. When this happens we can always find
             a germ of a hypoanalytic function at $p$ such that $\mathrm{d\,Z}(p)\neq 0$ such that
              $\{\mathrm{Z}=\mathrm{Z}(p)\}=\{p\}$. Thus being nowhere strictly pseudoconvex is necessary condition for a hypoanalytic structure verify the local Hartogs property.

\vskip0.5cm
We now adress the question of whether the $global$ Hartogs property holds for all pairs $(K,U)$ of compact sets $K\subset{U}$ where $U\subset\mathrm{M})$ is open. Hopefully we answer the question of Nacinovtch and Hill about
the example of the CR-structure on the hypersurface
 $|z_{1}|^{2}+|z_{2}|^{2}-|z_{3}|^{2}=1$ in $\textbf{C}^{3}$ where the zero of
the restriction of $z_{3}$ becomes compact failing the $global$ Hartogs property but curiously holding the $local$ one.
Such hypersurface is actually a zero of a homogenous solution of a Tube structure globally defined by the map
$$\mathrm{Z}:\textbf{R}^{2}\times\textbf{R}^{2}\times\textbf{R}^{2}\times\textbf{R}^{2}\longrightarrow\textbf{C}\times\textbf{C}\times\textbf{C}\times\textbf{C}$$
where $\mathrm{Z}(v_{1},v_{2},v_{3},v_{4} )=(z_{1},z_{2},z_{3},x+\mathrm{i\,}\Phi(|v_{1}|,|v_{2}|,|v_{3}|)$ with
$\mathrm{d\,z}_{1}\wedge\mathrm{d\,z}_{2}
\wedge\mathrm{d\,z}_{3}
\wedge\mathrm{d\,x}+\mathrm{i\,}
\mathrm{d\,}\Phi\neq 0$
on $\textbf{R}^{2}\times\textbf{R}^{2}
\times\textbf{R}^{2}\times\textbf{R}^{2}$. Finally it is taken
$\Phi(\xi_{1},\xi_{2},\xi_{3})=\xi^{2}_{1}+\xi^{2}_{2}-\xi^{2}_{3}$ and the zeroes of $\{x+\mathrm{i\,}\Phi(|v_{1}|,|v_{2}|,|v_{3}|)\}$
becomes CR-substructures which are actually also Tubes. By means of a right biholomorphism we find
$\mathrm{Z}(v_{1},v_{2},v_{3},v_{4} )=(z_{1},z_{2},z_{3},x+\mathrm{i\,}\Phi({\mathrm{Im\,}}z_{1},\mathrm{Im\,}z_{2},\mathrm{Im\,}z_{3}))$ thus a tube according the Definition VI.9.2 in Treves([Tr1]) and the embedded CR-submanifolds
$\{x+\mathrm{i\,}\Phi(|v_{1}|,|v_{2}|,|v_{3}|)=constant\}$  are Tubes with the restriction of $\mathrm{Z}_{0}=(z_{1},z_{2},z_{3})$ as a global integral.
After a unitary linear right composition the intersection zeroes
$\{x+\mathrm{i\,}\Phi({\mathrm{Im\,}}z_{1},\mathrm{Im\,}z_{2},\mathrm{Im\,}z_{3})=z_{3}=a+\mathrm{i\,}b\}$ have the expression
$$\big\{\textbf{R}^{2}\times\{(t_{1},t_{2})\in\textbf{R}^{2}:{\mathrm{Im\,}}^{2}z_{1}+\mathrm{Im\,}^{2}z_{2}=c+b^{2}\}\times\{a+\mathrm{i\,}b\}\times
\{a+\mathrm{i\,}c\big\}$$
and it is empty if $c<-b^{2}$. When $c=b^{2}$ it is the plane
$$\big\{\textbf{R}^{2}\times\{(0,0)\}\times\{a+\mathrm{i\,}b\}\times
\{a+\mathrm{i\,}c\big\}$$ becoming
homeomorphic to
$\textbf{R}^{2}\times\textbf{S}$ for all $c>-b^{2}$. It happens that
the function  $\mathrm{Re\,}[\mathrm{i\,}z_{3}+\kappa (z^{2}_{1}+z^{2}_{2})]=(1-\kappa)[\mathrm{Im\,}^{2}z_{1}+\mathrm{Im\,}^{2}z_{2}]+
\kappa[\mathrm{Re\,}^{2}z_{1}+\mathrm{Re\,}^{2}z_{2}]$ has for  $0<\kappa<1$ compact zeros homeomorphic to $\textbf{S}^{3}\subset
\textbf{R}^{2}\times\textbf{R}^{2}\times\{a+\mathrm{i\,}b\}\times\{a+\mathrm{i\,}\Phi({\mathrm{Im\,}}z_{1},\mathrm{Im\,}z_{2},b)=a+\mathrm{i\,}c\}$ and
noncompact zeros for $\kappa\geq 1$.

Let us now denote by $\mathcal{C}_{\mathrm{Z}}(p)$ the closed set $\{\Re\,{\mathrm{Z}\leq\Re\,\mathrm{Z}}(p)\}$ for an arbitrary hypoanalytic function $\mathrm{Z}$.
We will say that
 a Tube structure $(\mathrm{M}, \mathcal{L},\mathcal{G})$ verifies
 the $global$ Hartogs condition $(\mathtt{H})$ if:

$\bullet_{4}$ the Lie algebra $\mathcal{G}$ admits at least $\mathrm{codim\,}\mathcal{L}$
globally solvable generators,

$\bullet_{5}$ $\mathcal{C}_{\mathrm{Z}}(p)\subset\mathrm{M}$ does not have  compact components for all
hypoanalytic

 function $\mathrm{Z}$ and $p\in\mathrm{M}$,

$\bullet_{6}$ $\mathcal{O}_{\mathcal{L}}(p)\subset\mathrm{M}$ is never compact for all
$p\in\textrm{M}$.

The example of Hill $\&$ Nacinovitch ([HN])  will show that condition $\bullet_{4}$ is necessary. The condition $\bullet_{5}$ is obviously needed for
global Hartogs property holds.
Otherwise any open set containing
 one compact component would fail the Hartogs property. Finally  we may consider the quotient space $\mathcal{O}_{\mathcal{L}}$ defined by the equivalent relation $\sim$, where $p\sim q$ in $\textrm {M}$ if and only if $p,q\in \mathcal{O}_{\mathcal{L}}(p)$. Then every  $u\in\textrm{C}(\mathcal{O}_{\mathcal{L}})$ can be lifted to a
 function in $\textrm{C}(\textrm{M})$ which is a weak solution of ${\mathcal{L}}$ and $(u-u(p))^{-1}$ fails the global hartogs property showing that    $\bullet_{6}$ is also necessary.

  We can now state that

$\mathbf {Theorem. B.}$ Let  $\textrm{M}$ be simply connected and $(\mathrm{M}, \mathcal{L}, \mathcal{G})$ be a Tube structure.
Then $(\mathrm{M}, \mathcal{L}, \mathcal{G})$ verify $(\texttt{H})$ if and only if $global$
Hartogs property holds.

{\bf Remark.}
 This gives a explanation for the embedded example of
Hill $\&$ Nacinovitch (see [HN]) which gives an example of a Tube structure
which verify the local Hartogs phenomena but not the global one.
The Tube structure in question is defined by  given by the map $$\mathrm{Z}:\textbf{C}^{3}\times\textbf{R}\longrightarrow\textbf{C}^{4}$$ defined by $\mathrm{Z}(z,y)=\big(z,y+\mathrm{i\,}(|z^{2}_{1}|+|z^{2}_{2}|-|z^{2}_{3}|)\big).$
By means of an biholomorphism in $\textbf{C}^{4}$ we may rewrite
$\mathrm{Z}:\textbf{C}^{3}\times\textbf{R}\longrightarrow\textbf{C}^{4}$ as
$$\mathrm{Z}(z,x)=(x_{1}+\mathrm{i\,}t_{1},x_{2}+\mathrm{i\,}t_{2},x_{3}
+\mathrm{i\,}t_{3}, x+\mathrm{i\,}(t_{1}^{2}+t_{2}^{2}-t_{3}^{2})).$$
In this case there exists only one orbit and every the Hartogs global phenomena holds. On the other hand the zero $\mathcal{C}_{\mathrm{Z}}(p)$ with
$\mathrm{Z}(p)=\mathrm{i\,}$ is a hypoanalytic submanifold which happens to be a globally integrable.
  The global integral in question is the restriction of
  $(x_{1}+\mathrm{i\,}t_{1},x_{2}+\mathrm{i\,}t_{2},x_{3}
+\mathrm{i\,}t_{3})$ to the zero $$ \mathcal{C}_{\mathrm{Z}}(p)=\{x+\mathrm{i\,}
(t_{1}^{2}+t_{2}^{2}-t_{3}^{2}))=\mathrm{i\,}\}.$$
 Thus it is a Tube which enjoy the the local Hartogs phenomena but not the global one. It happens that  $(x_{3}
+\mathrm{i\,}t_{3})^{-1}$ is well defined in $\mathcal{C}_{\mathrm{Z}}(p)$ with
$\mathrm{Z}(p)=\mathrm{i\,}$ except by its with intersection with $x_{3}
+\mathrm{i\,}t_{3}=0$. The latter intersection is a set homeomorphic to the cylinder $\textbf{R}^{2}\times\textbf{S}^{1}$ and one can check that the function
$$\mathrm{Z}_{\kappa}=x_{3}+\mathrm{i\,}t_{3}+\kappa\,((x_{1}+\mathrm{i\,}t_{1})^{2}+(x_{2}+
\mathrm{i\,}t_{2})^{2})$$ for some small $k<1$
has a compact zero inside a Torus contained in $\textbf{R}^{2}\times\textbf{S}^{1}\subset\textbf{R}^{2}
\times\mathrm{i\,}\textbf{R}^{2}\simeq\textbf{C}^{2}$ violating the condition $\bullet_{5}$
 in Theorem B. The characterization given in [HT] for the global Hartogs phenomena in
 here stands for $\bullet_{5}$ one of the global condition in $(\mathtt{H})$. Thus the  Theorem B
 is a generalization of the result presented there.

\centerline{\bf 2.Proofs of Theorem A and B}

{\bf{Proof of Theorem A.}}
        It follows from the main result in [HT] that a Tube structure $(\mathrm{M},\mathcal{L},\mathcal{G},)$ enjoys the $local$
             Hartogs property if and only if the germ
             $\mathcal{C}_{\mathrm{Z}}(p)$  of a hypoanalytic function
             $\mathrm{Z}$ with $\mathrm{d\,Z}(p)\neq 0$ at $p$ is connected and  equivalent to nowhere strictly pseudoconvexity for Tubes structures $(\mathrm{M},\textsc{a},\mathcal{L})$.
             Observe that  $\textbf{C}\otimes
             \mathrm{T}_{p}\mathrm{M}=\mathcal{L}_{p}+
             \textbf{C}\otimes\mathcal{A}_{p}\subset\textbf{C}\otimes\mathrm{T}_{p}\mathcal{O}_{\mathcal{L}}+
             \textbf{C}\otimes\mathcal{A}_{p}$
             and consequently $\textbf{C}\otimes\mathrm{N}^{\ast}\mathcal{O}_{\mathcal{L}}(p)
             \subset\mathcal{L}^{\perp}$.
              If the  local Hartogs phenomena occurs for this hypoanalytic  structure then
the germ  $\mathcal{C}_{\mathrm{Z}}(p)$ for a hypoanalytic function with
 $\mathrm{d\,Z}(p)\neq 0$ must be connected. Otherwise
it will display some compact component or a denumerable set of components.
In the first case  $\mathrm{Z}^{-1}$ would fails the Hartogs
phenomena for some pair $(U,\mathcal{C}_{\mathrm{Z}}(p)\cap{U}\,)$ and in the second is void, otherwise
$\mathrm{d\,Z}(p)=0$.
By means of a complex linear transformation one may  assume that for a hypoanalytic chart $(\mathrm{Z}_{\alpha},U)$ with $p\in{U}$ that $\mathrm{Z}_{\alpha}(p)=0$ with
$\mathrm{d\,}\mathrm{Z}_{\alpha}(p)=\mathrm{I}$. Then
 $\mathrm{N}_{p}^{\ast}\mathcal{O}_{\mathcal{L}}(p)$ will necessarily have a basis among the differentials $\{\mathrm{d\,Re\,Z}_{1,\alpha}(p),...,\mathrm{d\,Re\,Z}_{m,\alpha}(p)\}$.
 This implies with $\mathrm{Z}_{\alpha}^{2}=
\mathrm{Z}_{1,\alpha}^{2}+\cdot\cdot\cdot+\mathrm{Z}_{m,\alpha}^{2}$ and large  $\kappa$ that
the germ
$\mathcal{C}_{\mathrm{Z}+\kappa\,\mathrm{Z}_{\alpha}^{2}}(p)
\cap\mathcal{O}_{\mathcal{L}}(p)$ must be connected if $\mathcal{C}_{\mathrm{Z}}(p)$ is.
 Consequently for Tubes a necessary and sufficient condition the validity of local Hartogs phenomena is translated  on germs  $\mathcal{F}_{\mathrm{Z}}(p)=
\mathcal{C}_{\mathrm{Z}}(p)\cap\mathcal{O}_{\mathcal{L}}(p)$ by the condition;
\centerline{{\it  $\mathcal{F}_{\mathrm{Z}}(p)$ is connected for all hypoanalytic germs $\mathrm{Z}$ with
 $\mathrm{Ker\,d\,Z}(p)= \{0\}$} \,\,\,\,\,\, ($\mathtt{P}$).}

\underline{\underline{}} {\textbf{Proof of Theorem B.(1st version via Ehrenpreis argument )}}

We will prove that the first cohomological group of the complex induced by $\mathcal{L}$ is trivial reviving the original idea of Ehrenpreis [Eh] and giving a stronger version of Theorem.B.
We  select $m=\textrm{codim}\mathcal{L}$ globally integrable vector fields from $\mathcal{G}$ and assume without loss of generality that $\textrm{dim}\mathcal{G}=m$.  It follows that there exist smooth manifold $\textrm{N}$ such that
$\textrm{M}=\textbf{R}^{m}\times\textrm{N}$.  For $m=1$ it is the statment of Theorem 6.4.2 (f) in [DH]. When $m$ is bigger than one we proceed by induction
taking advantage of the commutative property of the fields.
   As a consequence we get an open projection $\Pi_{\textsc{a}}: \textrm{M}\rightarrow\textrm{N}$ having as fibers the the $m-$dimensional submanifolds $\mathrm{A}\subset\mathrm{M}$ verifying $\mathrm{T}_{p}\,\mathrm{A}=\mathcal{A}_{p}$ if $p\in\mathrm{A}$, that is $\mathrm{A}=\mathbf{R}^{m}\times\{\Pi_{\textsc{a}}(p)\}$.
 Now follows from the characterization of
tubes structures found in VI.8 Partial Local Group Structures([Tr1]) that one can construct a hypoanalytic  atlas $(U_{\alpha}, \mathrm{Z}_{\alpha})$ such that $\mathrm{Z}_{\alpha}(x,t)=x+\textrm{i}\Phi(t)$ where the first coordinates $x$ are first integrals of the chosen $m$ globally integrable vector fields in $\mathcal{G}$.
Since  a pair of  hypoanalytic charts $(U_{\alpha}, \mathrm{Z}_{\alpha}),(U_{\alpha^{\prime}}, \mathrm{Z}_{\alpha^{\prime}})$ changes by a biholomorphism and they have identical real parts in $U_{\alpha}\cap U_{\alpha^{\prime} }$ they must agree there.
It follows that $\mathcal{L}$ has a global integral $\mathrm{Z}$ and
  the topological space $\mathrm{M}/\sim_{\mathcal{F}}$, where $\sim_{\mathcal{F}}$ is the equivalence of being in a same fiber of $\mathcal{L}$, is globally defined.
Let $\mathrm{Z}=(\mathrm{Z}^{1},...,\mathrm{Z}^{m})$ be a global integral for $(\mathrm{M}, \mathcal{L}, \textsc{a} )$,  that is a map from $\mathrm{M}\longrightarrow\textbf{C}^{k}$ with
 $k=\mathrm{codim\,}\mathcal{L}$. Let $\omega=\sum_{j=1}^{n}f_{j}\textrm{d}t_{j}$ be a smooth closed class in the first cohomological group with compact support induced by the differential complex associated to $\mathcal{L}$. We mean that $$\textrm{d}\omega\wedge\textrm{d}\textrm{Z}=0$$ where
$\textrm{d}\textrm{Z}=\textrm{d}\textrm{Z}_{1}\wedge...\wedge\textrm{d}\textrm{Z}_{m}$.
 The same steps in [HT] by performing Fourier transform of $\omega\wedge\textrm{d}\textrm{Z}$ in the linear fibers $\{t\}\times\textbf{R}^{m}$ to find $\hat{\omega}\in\wedge^{1}\textrm{T}^{\ast}(\textrm{N})$ such that
 $$\textrm{d}_{t}\textrm{e}^{\Phi\cdot\xi}\hat{\omega}=0\,\,\,\hbox{for all}\,\,\,\xi\in{\textbf{R}^{m}}^{\ast}.$$
 Since $\textrm{M}$ is simply connected so it is $\textrm{N}$ which enables us to define
 $v(\xi,t)$ by $$\textrm{d}_{t}v(\xi,t)=\textrm{e}^{\Phi\cdot\xi}\hat{\omega}, \,\,\,\,\,v(\xi,t_{0})=0$$
 where $t_{0}\in \textrm{N}\setminus \Pi_{\textrm{A}}(\textrm{supp\,}{\omega})$. Now set $\hat{u}(\xi,t)=\textrm{e}^{-\Phi\cdot\xi}{v}(\xi,t)$ which vanishes outside $\Pi_{\textrm{A}}(\textrm{supp\,}{\omega})$. It remains to prove that $\hat{u}$ is indeed the fiber Fourier transform of a function
 $u\in\textrm{C}^{\infty}_{c}(\textrm{M})$ to finishes the proof.
It follows from $\bullet_{5}$ that the sublevels
 $$\mathcal{C}_{-\mathrm{i Z}}(0,t)=\{s\in\textrm{N}:\Phi(s)\cdot\xi\leq \Phi(t)\cdot\xi\}$$ does not have compact components. Since it is a closed set it implies that $\textrm{N}$ can not be compact. We now cover $\textrm{N}$ by charts
 $(\chi_{\beta}, W_{\beta})$ associated to the maximal atlas of $\textrm{N}$ such that each one maps $W_{\beta}$ onto   $Q_{0}=[0,1]^{n}$. Also we may assume that $\{W_{\beta}\}$ is a locally finite covering. Now we consider a subdivision of $Q_{0}$ in  $2^{nk}$ cubes $Q_{k}$ of side length $2^{-k}$. Then any polygonal line inside $Q_{0}$ which intercepts each division cube in a unique line segment will have a length bounded by $\sqrt{n}2^{-k}2^{nk}=\sqrt{n}2^{(n-1)k}$. In particular the image of a such polygonal line  by $\chi^{-1}_{\beta}$ into $W_{\beta}$ will have length bounded by $C_{\beta}\sqrt{n}2^{(n-1)k}$ for some metric in $\textrm{N}$ which is equivalent to the euclidian metric of $Q_{0}$ via any $\chi_{\beta}$.
 We now consider only  cubes $Q_{k}$ such that $\chi_{\beta}(Q_{k})$ meets the $connected$ component of $\mathcal{C}_{-\mathrm{i Z}}(0,t)$ which contains $t$  for some $\beta$.
 It entails that $\cup_{\beta}\chi_{\beta}(Q_{k})\supset\mathcal{C}_{-\mathrm{i Z}}(0,t)$ is a connected set and we can find a curve differentiable by parts $\gamma$ linking $t$ to an arbitrary point in $\cup_{\beta}\chi_{\beta}(Q_{k})$
 such that $\chi_{\beta}(\gamma)$ is a polygonal curve in $[0,1]^{n}$ which meets any $Q_{k}$ in a line segment for all
 $\beta$. Now every $s\in\gamma$ is at a distance (for the chosen metric) comparable with $\sqrt{n}2^{-k}$
 from the component of  $\mathcal{C}_{-\mathrm{i Z}}(0,t)$ which contains $t$. Let $t^{\prime}$ a point of the component within this range and apply the mean value theorem to obtain
 $$|\Phi(s)-\Phi(t^{\prime})\cdot\xi|\leq \sup_{t^{\prime}\in{Q}_{k}}|\nabla\Phi(t^{\prime})|\sqrt{n}2^{-n}|\xi|.$$
 It is also true that $(\Phi(s)-\Phi(s^{\prime}))\cdot\xi\leq 0$ for every $\xi\in{\{t\}\times\textbf{R}^{m}}^{\ast}$
 with $t\in\textrm{N}$ and we can choose $\gamma$ such that
$\partial{Q_{k}}\cap\chi_{\beta}(W_{\beta}\cap{\gamma})$ oriented set $\{t_{0}, t_{1} \}$ obeying $|t_{0}-t|$ and $|t_{1}-t|$ are minimum and maximum of $|s-t|$ with $s\in \gamma\cap{Q_{k}}$. Now we may estimate $\hat{u}$ as $|\hat{u}(\xi,t)|=$
$$\bigg|\int_{\gamma}\textrm{e}^{(\Phi(s)-\Phi(t))\cdot\xi}\,\,\hat{\omega}\bigg|\leq\bigg|\int_{\gamma}
\textrm{e}^{(\Phi(t^{\prime})-\Phi(t))\cdot\xi}\,\hat{\omega}\bigg|\leq\sqrt{n}2^{(n-1)k}\textrm{e}^{C|\xi|2^{-k}}\sup|\hat{\omega}|$$
where the supreme of $|\hat{\omega}|$ is uniformly bounded in $\Pi_{\textrm{A}}(\textrm{supp\,}{\omega})$ by multiples of arbitrary
powers of $(1+|\xi|)^{-1}$.
Choosing $2^{-k}$ comparable with $(1+|\xi|)^{-1}$ we may find constants such that
$$|\hat{u}(\xi,t)|\leq C_{l}(1+|\xi|)^{-l}\,\,\hbox{for}\,\, \,\,t\in\Pi_{\textrm{A}}(\textrm{supp\,}{\omega})\,,\xi\in\textbf{R}^{m}\,,l\in{\textbf{N}}.$$
This happens because $\hat{u}(\xi,t)$ is uniformly bounded in the Schwartz space $\mathcal{S}(\textbf{R}^{m})$ for
every $t\in\textrm{N}$. It entails that $u(x,t)$, the Fourier inverse transform of $u(\xi,t)$ is indeed a function
in $\textrm{C}^{\infty}(\textrm{M})$. Compactness of $\textrm{supp}u$ follows from a theorem of propagations of zeroes of solutions for the sections of $\mathcal{L}$. It states that solutions which vanishes in a neighborhood of a point
$p\in\mathcal{O_{\mathcal{L}}}(p)$ must vanishes in all orbit(see Theorem 1.1 in [HP]). In our case we consider
the structure $(\textrm{M}\setminus\textrm{supp}u, \mathcal{L},\textrm{A})$ to apply the cited theorem. Uniqueness of the solution $u$ follows in a similar argument.

\underline{\underline{}} {\textbf{Proof of Theorem B.(2nd version via Arens-Royden theorem)}}
We  select $m=\textrm{codim}\mathcal{L}$ globally integrable vector fields from $\mathcal{G}$ and assume without loss of generality that $\textrm{dim}\mathcal{G}=m$.  It follows that there exist smooth manifold $\textrm{N}$ such that
$\textrm{M}=\textbf{R}^{m}\times\textrm{N}$.  For $m=1$ it is the statment of Theorem 6.4.2 (f) in [DH]. When $m$ is bigger than one we proceed by induction
taking advantage of the commutative property of the fields.
   As a consequence we get an open projection $\Pi_{\textsc{a}}: \textrm{M}\rightarrow\textrm{N}$ having as fibers the the $m-$dimensional submanifolds $\mathrm{A}\subset\mathrm{M}$ verifying $\mathrm{T}_{p}\,\mathrm{A}=\mathcal{A}_{p}$ if $p\in\mathrm{A}$, that is $\mathrm{A}=\mathbf{R}^{m}\times\{\Pi_{\textsc{a}}(p)\}$.
 Now follows from the characterization of
tubes structures found in VI.8 Partial Local Group Structures([Tr1]) that one can construct a hypoanalytic  atlas $(U_{\alpha}, \mathrm{Z}_{\alpha})$ such that $\mathrm{Z}_{\alpha}(x,t)=x+\textrm{i}\Phi(t)$ where the first coordinates $x$ are first integrals of the chosen $m$ globally integrable vector fields in $\mathcal{G}$.
Since  a pair of  hypoanalytic charts $(U_{\alpha}, \mathrm{Z}_{\alpha}),(U_{\alpha^{\prime}}, \mathrm{Z}_{\alpha^{\prime}})$ changes by a biholomorphism and they have identical real parts in $U_{\alpha}\cap U_{\alpha^{\prime} }$ they must agree there.
It follows that $\mathcal{L}$ has a global integral $\mathrm{Z}$ and
  the topological space $\mathrm{M}/\sim_{\mathcal{F}}$, where $\sim_{\mathcal{F}}$ is the equivalence of being in a same fiber of $\mathcal{L}$, is globally defined.

Let
   $\mathrm{Z}$ be a global integral for $(\mathrm{M}, \mathcal{L}, \textsc{a} )$,  that is a map from $\mathrm{M}\longrightarrow\textbf{C}^{k}$ with
 $k=\mathrm{codim\,}\mathcal{L}$.
 Now, under the hypothesis $(\mathtt{P})$ the closed set $\mathcal{Z}(p)=
              \{\mathrm{Z}=\mathrm{Z}(p)\}$ is
              locally connected and the map $\Pi
 _{\mathcal{L}}:\mathrm{M}/\sim_{\mathcal{L}}\,\,
 \hookrightarrow\mathrm{Z(\mathrm{M})}$ is relatively open and locally injective.
               Here
$\sim_{\mathcal{L}}$ represents the equivalence relation of two points of  $\mathrm{M}$
being in a same component of the closed  set
         $\mathcal{Z}(p)$, that is
 the set $\mathcal{Z}(p)$ agree locally with $\mathcal{F}_{p}$ implying that $\mathcal{Z}_{\mathrm{Z}}(p)=
                \cup_{p\in\mathcal{Z}(p)}\mathcal{F}_{p}$, thus invariantly defined.  We call $\mathrm{M}/\sim_{\mathcal{L}}$ the reduced manifold by $\mathcal{L}$ which makes any automorphism commute with a homeomorphism of $\mathrm{M}/\sim_{\mathcal{L}}$ via the canonical projection. Such subgroup of homeomorphisms is a hypoanalytic invariant since the germ of the fiber $\mathcal{F}({p})$ propagates trough $\mathcal{Z}(p)$. Thus we may say that
                the fiber of $\mathcal{L}$ is globally defined and
                  $\mathrm{M}/\sim_{\mathcal{L}}$ is invariant under automorphism of the structure. We mean by global diffeomorphisms of $\mathrm{M}$ which leaves $\mathcal{L}$ invariant in the sense that its differential is an automorphism of $\mathcal{L}_{p}$ for
                every $p\in\mathrm{M}$. We say that an open subset $U\subset\mathrm{M}$ is a domain for $\mathcal{L}$ if the canonical projection
                $\Pi_{\mathcal{L}}:\mathrm{M}\rightarrow\mathrm{M}
                /\sim_{\mathcal{L}}$ is injective in $U$. If the intersection $\mathcal{C}_{\mathrm{Z}}(p)
                \cap\mathcal{O}_{\mathcal{L}}(p)$
              is relatively open in $\mathcal{O}_{\mathcal{L}}(p)$ then by uniqueness(see [Tr1]) $\mathcal{O}_{\mathcal{L}}(p)\subset \mathcal{C}_{\mathrm{Z}}(p)$ and the germ propagates into the orbit  $\mathcal{O}_{\mathcal{L}}(p)$. Despite the
discreteness of fibers of the canonical projection $\Pi_{\sim_{\mathcal{L}}}$, one can not expect that
               $\mathrm{M}/\sim_{\mathcal{L}}
               $ evenly covers $\mathrm{Z(\mathrm{M})}$ and in this way not necessarily a covering space.
              Now, for any compact subset $K\subset\mathrm{M}/\sim_{\mathcal{L}}$ we
we consider its $\mathcal{L}-$convex envelope
 $\widehat{K}\subset\mathrm{M}/\sim_{\mathcal{L}}$ with respect to the finitely generated Banach algebra $\mathcal{A}_{\mathcal{L}}(K)$ of continuous functions $u$ of ${K}$ which are uniform limits in $K$ of polynomials in $\mathrm{Z}$. Such continuous are of course also defined in
$\widehat{K}$ the polynomial convex envelope of $K$ and $\mathrm{Z}(\widehat{K})=\widehat{\mathrm{Z}({K})}$.  Thus
$\mathrm{Z}(\widehat{K})$ indeed agree with the maximal ideal space of the algebra $\mathcal{A}_{\mathcal{L}}(K)$(see Theorem 3.1.15 in ([Ho])). If the first Cech cohomology group of $\mathrm{H}^{1}({K}/\sim_{\mathcal{L}})$ is not trivial it follows from the Arens-Royden theorem (see [Ar], [Ro]) that we can find a  hypoanalytic polinomial $\mathrm{Z}_{0}$ such that $\mathrm{d\,Z}_{0}/\mathrm{Z}_{0}\neq 0$ is well defined in $\mathrm{M}$ and a representant  of a non trivial class in
 $\mathrm{H}_{\mathrm{d\,}}^{1}(\mathrm{M})$ (the first cohomological DeRham group of the complex defined by the exterior derivative $\mathrm{d\,}$). If the DeRham cohomology group $\mathrm{H}_{\mathrm{d\,}}^{1}(U\setminus K)$ is trivial then
 $K\subset \mathrm{M}$ is an irremovable singularity of the ring $\mathcal{A}(\mathrm{M})$ because in this case
 $\mathrm{Log\,} \mathrm{Z}_{0}$ will be a hypoanalytic function which is defined in $\mathrm{M}\setminus K$ which cannot be extended for all
$\mathrm{M}$ failing the Hartogs phenomena for the pair $(\mathrm{M},K)$. On the other hand it follows from Poincar\'e duality that
$$\lim_{K\subset\subset\mathrm{M}}
\mathrm{H}_{\mathrm{d\,}}^{p}(\mathrm{M}, \mathrm{M}\setminus{K})
\simeq\mathrm{H}_{\mathrm{d\,}c}^{p}(\mathrm{M})
\simeq\mathrm{{H}}_{\mathrm{d\,}}^{\mathrm{dim \,M}-p}(\mathrm{M}).$$
Since in paracompact differentiable manifolds Cech, singular and De Rahm cohomology agree and $\mathrm{M}/\sim_{\mathcal{L}}$ inherits from the manifold $\mathrm{M}$ a CW-complex structure. It follows that
 the  Cech  and singular cohomology of
$\mathrm{M}/\sim_{\mathcal{L}}$ are well defined and agree.
If $\sim_{\mathcal{L}}$ is $proper$ then there exist
 a natural injection $$\mathrm{
{H}}_{c}^{p}(\sim_{\mathcal{L}}):
\mathrm{ {H}}_{c}^{p}(\mathrm{M}/\sim_{\mathcal{L}})\hookrightarrow
\mathrm{H}^{p}_{c}(\mathrm{M})\simeq
\check\mathrm{H}_{c}^{p}(\mathrm{M})$$ given by the singular cohomology functor. With $m+n=\mathrm{dim\,}\mathrm{M}$ we have
$$
\mathrm{H}_{c}^{m+1}(\mathrm{M}/\sim_{\mathcal{L}})
\simeq\mathrm{H}_{\mathrm{n-1}}(\mathrm{M}/\sim_{\mathcal{L}})$$

where $n=\mathrm{dim\,}\mathcal{L}$ by Poincar\'{e} duality. Now we have direct
decomposition  $$\mathrm{H}_{\mathrm{d\,}c}^{m+1}(\mathrm{M})\simeq
\mathrm{H}_{c}^{m+1}(\sim_{\mathcal{L}})
[\mathrm{H}_{c}^{m+1}
(\mathrm{M}/\sim_{\mathcal{L}})]
\oplus\mathrm{Ker\,}\wedge\Omega$$ where
$\Omega=\mathrm{d\,}\zeta$ is a exact nonvanishing section of $\wedge^{m}\mathcal{L}^{\perp}$ and the

$$\wedge\Omega:\wedge^{1}\mathrm{T}^{\ast}(\mathrm{M})
{\longrightarrow}\wedge^{m+1}
\mathrm{T}^{\ast}(\mathrm{M})$$ is defined for $\omega\in\wedge_{c}^{1}\mathrm{T}^{\ast}(\mathrm{M})
$ by $\omega\wedge\Omega$ verifies
$\Omega\wedge\mathrm{d}=\mathrm{d}\wedge\Omega$ and induces homomorphism
$\wedge\Omega:\mathrm{H}_{\mathrm{d\,}c}^{1}(\mathrm{M})\rightarrow
\mathrm{H}_{\mathrm{d\,}c}^{m+1}(\mathrm{M})$. We can represent $\mathrm{H}^{1}_{{\mathrm{d}_{\mathcal{L}}c}}
(\mathrm{M})$ (where ${\mathrm{d}_{\mathcal{L}}}$ is the exterior derivative induced by $\mathrm{d}$ in the sections of
$\textbf{C}\otimes\mathrm{T}^{\ast}\mathrm{M}/\mathcal{L}^{\perp}$) as the kernel of the map
$\omega\mapsto\omega\wedge\Omega$ in $\mathrm{H}_{\mathrm{d\,}c}^{1}(\mathrm{M})$.
 Thus $\omega\wedge\Omega$ represent a class in
$\mathrm{H}_{\mathrm{d\,}c}^{m+1}(\mathrm{M})
$ if $\omega$ is represents a class in $\mathrm{H}^{1}_{{\mathrm{d}_{\mathcal{L}}c}}
(\mathrm{M})$. In this setting
the Hartogs phenomena holds if and only if  for all
$\omega\in{\mathrm{H}_{\mathrm{d}_{\mathcal{L}\,}c}^{1}}
(\mathrm{M})$ there exist $u\in \mathrm{C}_{c}^{\infty}(\mathrm{M})$ such that
$$\mathrm{d}u\wedge \Omega=\omega\wedge \Omega \eqno(\ast)$$

Solvability of $(\ast)$ assures the triviality of the intersection
$\mathrm{H}_{\mathrm{d\,}c}^{m+1}(\mathrm{M})\cap
\mathrm{H}^{
1}_{{\mathrm{d}_{\mathcal{L}}}}(\mathrm{M})=
\mathrm{H}_{\mathrm{d}_{\mathcal{L} \,}c}^{1}(\mathrm{M})$  which in turn must represent some subgroup of the de Rham group $\mathrm{H}_{\mathrm{d\,}}^{n-1}(\mathrm{M})$ via Poincar\'{e} duality. The existence of  a Lie algebra $\textsc{a}$ oriented by $\Omega$ allows one to decompose
$\mathcal{L}\subset\textbf{C}\otimes\textsc{a}
\oplus\mathrm{T}\mathrm{B}$ where
$\mathrm{B}=\Pi_{\textsc{a}}(\mathrm{M})$ is a real $n-$dimensional manifold obtained by identifying  the fibers of $\Pi_{\textsc{a}}$ to points in $\mathrm{M}$.
It follows that  every $real$ section of $\mathrm{T}\mathrm{B}$ has a unique lifting to $\mathcal{L}$.
This enables us to define the connection $$\nabla_{T}L (p)=\mathrm{T}(L)(T(p))-T_{h}(L(p))\in\textbf{C}\otimes\textsc{a}_{p}$$
where $\mathrm{T}\Pi_{\textsc{a}}(T_{h})=T$ at $p$ when the fibers $\mathcal{A}_{p}$ of $\Pi_{\textsc{a}}$ have a affine linear structure,
and this is always the case for a open covering $U_{\alpha}$ of
$\mathrm{M}$ such that $\textsc{a}$ admits $m-1$ globally solvable
generators in
$\Pi_{\textsc{a}}^{-1}(\Pi_{\textsc{a}}(U_{\alpha}))$, turning  $\mathrm{M}$ into a real vector bundle by defining local charts $\Pi_{\textsc{a}}^{-1}(\Pi_{\textsc{a}}(U_{\alpha}))\simeq
\textbf{R}^{m}\times\Pi_{\textsc{a}}(U_{\alpha})$.
Assume that
$\mathrm{H}^{m+1}_{\mathrm{d\,}}(\mathrm{M})\simeq{\mathrm{H}^{m+1}}{(\mathrm{M})}$ verifies $\mathrm{H}^{m+1}_{\mathrm{d\,}}(\mathrm{M})
=\{0\}$ which means that any section $\omega\wedge\Omega$ is automatically exact if it represents  a class in $\mathrm{H}^{1}_{{\mathrm{d}_{\mathcal{L}}}}(\mathrm{M})$. Then we can find a section $e\,\Omega+ \lambda $ of $\wedge^{m}\mathrm{T}^{\ast}(\mathrm{M})$ such that $\mathrm{d\,}(e\,\Omega+ \lambda)=\omega\wedge\Omega$.
It follows from the Stoke's
Theorem that for for rectifiable $m+1-$ rectifiable chain of form  $\sigma=\Pi_{\textsc{a}}^{-1}(\Pi_{\textsc{a}}(\sigma))
=\textbf{R}^{m}\times{\Pi_\textsc{a}}(\sigma)$ that
$$\int_{\partial\sigma}e\,\Omega=
\int_{\partial\sigma}(e\,\Omega+ \lambda)=\int_{\sigma}\omega\wedge\mathrm{d\,}\Omega=
\int_{t\in\Pi_{\textsc{a}}
(\sigma)}\int_{\Pi^{-1}(t)}
\omega\wedge\Omega=
\int_{\Pi_{\textsc{a}}
(\sigma)}\int_{\textbf{R}^{m}}
\omega\wedge\Omega
$$
for all  $\omega\in\mathrm{H}_{\mathrm{d\,}c}^{m+1}(\mathrm{M})$.
 In particular $\sigma$ is invariant by the $\textsc{a}-$flow and the left side is finite if $\omega$ has compact support.
 Thus if $\sigma$
is a $m+1-$ rectifiable chain with  boundary $\partial\sigma$ and
$\Omega_{|\sigma}\neq 0$ then  locally $\Pi_\textsc{g}(\sigma)$ is a $1-$rectifiable.  If we choose $\sigma$ such that $\Pi_{\textsc{a}}
(\partial\sigma)=\{t\}$ then the left side above is a smooth function of  $t$ which vanishes outside $\Pi_{\textsc{a}}(\mathrm{supp\,}\omega)$. We finish the proof applying the Treves propagation of zeroes theorem as we did before.

\end{document}